\documentclass[12pt]{article}
\usepackage{amsmath}
\usepackage{latexsym,amssymb,eucal}
\newtheorem{thm}{Theorem}[section]
\newtheorem{lem}[thm]{Lemma}

\newtheorem{df}[thm]{Definition}
\newtheorem{exam}[thm]{Example}
\newcommand{\id}{\mathrm{id}}
\newcommand{\Ad}{\mathrm{Ad}\,}

\newcommand{\cM}{\mathcal{M}}
\newcommand{\cN}{\mathcal{N}}
\newcommand{\cB}{\mathcal{B}}
\newcommand{\cA}{\mathcal{A}}

\newcommand{\Aut}{\mathrm{Aut}}
\newcommand{\Sect}{\mathrm{Sect}}
\newcommand{\End}{\mathrm{End}}
\newcommand{\Ker}{\mathrm{Ker}}
\newcommand{\Mor}{\mathrm{Mor}}
\newcommand{\Int}{\mathrm{Int}}

\begin{document}

\title{A simple sufficient condition for \\ 
triviality of obstruction in  \\ 
the orbifold construction for subfactors}
\author{Toshihiko MASUDA\footnote{Supported by 
JSPS KAKENHI Grant Number 23540246.} \\
 Graduate School of Mathematics, Kyushu University \\
 744, Motooka, Nishi-ku, 
Fukuoka, 819-0395, JAPAN}
\date{}

\maketitle
\begin{abstract}
We present a simple sufficient condition for triviality of
obstruction in the orbifold construction. As an application, we can show the existence
of subfactors with principal graph $D_{2n}$ without full use of
Ocneanu's paragroup theory.
\end{abstract}

\section{Introduction}\label{sec:intro}
In the subfactor theory initiated by V. F. R. Jones \cite{J-ind}, one of
the fundamental construction of subfactors is the orbifold construction. 
It was introduced by Kawahigashi \cite{Kw-orbi} to show the
existence of subfactors whose principal graphs are Dynkin diagram
$D_{2n}$. Roughly speaking, the orbifold construction is to take a
``quotient'' by an internal symmetry of subfactors, which is realized by
taking a crossed product construction by an abelian group. 
The orbifold construction has been further studied by
\cite{EK-HeckeI}, \cite{Xu-orb}, \cite{Xu-flat}, \cite{Go-orbi},
\cite{Go-Loi}.

A typical example of the orbifold construction is the case of Jones
subfactors with principal graph $A_{2n-1}$ \cite{J-ind},
which possess certain $\mathbb{Z}/2\mathbb{Z}$-symmetry. 
By the orbifold construction, graph change occurs only for $A_{4n-3}$-subfactors, and
the orbifold construction produces subfactors with principal graph $D_{2n}$. This is because
some obstruction, which prevents graph change,  appears for
an $A_{4n-1}$-subfactor.
Thus the most important problem is to
determine the triviality of obstruction in the construction. 
In general, this problem requires complicated combinatoric computation of
connection.

In this paper, we present a simple sufficient condition for triviality of
obstruction appearing in the  orbifold construction. Namely, we show
that an obstruction vanishes if a tensor
category arising from a subfactor has a nice fusion rule.
As an application, we can show the existence
of subfactors with principal graph $D_{2n}$ without full use of
Ocneanu's paragroup theory, and the proof is easier than that of  \cite{Kw-orbi}.
 Our argument is inspired by the computation presented in
 \cite{Iz-CuntsI}, and we use sector technique for proof. 

This paper is organized as follows.
In \S \ref{sec:half}, we recall the definition of the Loi invariant
\cite{Loi-auto}, and construct a half braiding. This allows us to extend
an endomorphism to a crossed product factor, which is a special case of
$\alpha$-induction \cite{Xu-newbraided}, \cite{Boek-Evan-alphaI},
\cite{Boek-Evan-alphaII}, \cite{Boek-Evan-alphaIII},
\cite{BEK-alpha}. In \S \ref{sec:flat}, we 
show basic properties of this extension, and show the main theorem. Then
we present examples of application of our theorem.

\section{Half braidings}\label{sec:half}

We refer \cite{EK-book} for basic of subfactor theory, and 
\cite{Iz-fusion},  \cite{Iz-CuntsII} for basic of sector theory. In this
paper, we mainly treat type III factors. Here we recall some notations
of sectors.
For von Neumann algebras $\cA, \cB$, let $\Mor(\cA,\cB)$ be a set of all
unital continuous injective morphisms, and $\Sect(\cB,\cA)=\Mor(\cA,\cB)/\Int(\cB)$.
When $\cA=\cB$, we denote $\Mor(\cA,\cA)=\End(\cA)$, and
$\Sect(\cA,\cA)=\Sect(\cA)$.
For $\rho,\sigma\in \Mor(\cA,\cB)$, the space of intertwiners is defined
by $(\sigma,\rho):=\{a\in \cB\mid \rho(x)a=a\sigma(x), x\in \cA\}$. When
$\sigma$ is irreducible, i.e., $(\sigma,\sigma)=\sigma(\cA)'\cap
\cB=\mathbb{C}1$, $(\sigma,\rho)$ becomes a Hilbert space by an inner product
$\langle a,b\rangle 1 =b^*a$.

Let $\cN\subset \cM$  be an irreducible subfactor of type III  with finite index, and 
$\iota:\cN\hookrightarrow \cM$ an inclusion map. 
Then $\gamma=\iota \bar{\iota}$ is  a canonical
endomorphism for $\cN\subset \cM$.
Set 
\[
\cM\supset \iota(\cN) \supset \iota \bar{\iota}(\cM)\supset 
\iota\bar{\iota}\iota(\cN)\supset 
\cdots =
\cM\supset \cN_1\supset \cN_2\supset \cN_3\cdots. 
     \]

Set 
\begin{align*}
{}_\cM\Delta_{\cM}&=\{[\sigma]\in \Sect(\cM)\mid \sigma \prec (\iota\bar{\iota})^n, n\in
\mathbb{N}, \sigma \mbox{ is irreducible}\},  \\
{}_\cM\Delta_{\cN}&=\{[\sigma]\in \Sect(\cM,\cN)\mid \sigma \prec
(\iota\bar{\iota})^n\iota , n\in
\mathbb{N},\sigma \mbox{ is irreducible}\},  \\
{}_\cN\Delta_{\cM}&=\{[\sigma]\in \Sect(\cN,\cM)\mid \sigma \prec
 (\bar{\iota}\iota)^n\bar{\iota}, n\in 
\mathbb{N},\sigma \mbox{ is irreducible}\}, \\
{}_\cN\Delta_{\cN}&=\{[\sigma]\in \Sect(\cN)\mid \sigma \prec
(\bar{\iota}\iota)^n, n\in \mathbb{N},\sigma \mbox{ is irreducible}\}.   
\end{align*}
and $\Delta:={}_\cM\Delta_{\cM} \sqcup {}_\cM\Delta_{\cN} 
\sqcup {}_\cN\Delta_{\cM} 
\sqcup {}_\cN\Delta_{\cN}$.

We first recall the definition of the Loi invariant \cite{Loi-auto}. 
Fix isometries $R\in (\id, \iota\bar{\iota} )$ and 
$\bar{R}\in (\id, \bar{\iota}\iota )$ such that 
$\displaystyle{R^*\iota(\bar{R})=\bar{R}^*\bar{\iota}(R)=\frac{1}{d(\iota)}}$.
Let $\Aut(\cM,\cN)$ be a set of automorphisms of $\cM$ which preserve
$\cN$ globally. 
Take $\alpha\in \Aut(\cM,\cN)$.
Since  $\alpha$ preserves $\iota(\cN)$, 
we have $\alpha\iota=\iota \alpha$. Thus
$[\bar{\iota}][\alpha]=[\alpha][\bar{\iota}] \in \Sect(\cN,\cM)$, and we can
take $u\in U(\cN)$ such that $\alpha\bar{\iota}=\Ad (u)\circ \bar{\iota}\alpha$.
The choice of $u$ is not unique, but we can easily see
$\iota(u)^*\alpha(R)\in (\id,\iota\bar{\iota})$,
$u^*\alpha(\bar{R})\in (\id,\bar{\iota}\iota)$ and 
\[
 (u^*\alpha(\bar{R}))^*\bar{\iota}(\iota(u)\alpha(R))= 
(\iota(u)\alpha(R))^*\iota(u^*\alpha(\bar{R}))=\frac{1}{d(\iota)}.
\]
Thus we can choose a
unique $u_\alpha\in U(\cN)$ by $\iota(u_\alpha^*)\alpha(R)=R$, and 
$u^*_\alpha\alpha(\bar{R})=\bar{R}$. 
\begin{lem}\label{rem:1cocycle}
$(1)$ We have 
$\iota(u_{\alpha\beta})=\alpha\iota(u_\beta)\iota(u_\alpha)$ 
for  $\alpha,\beta\in \Aut(\cM,\cN)$. \\
$(2)$ We have $u_{\alpha}=v\bar{\iota}(v^*)$ for $\alpha=\Ad v$, $v\in U(\cN)$.
\end{lem}
\textbf{Proof.} The statement (1) follows from the uniqueness of
$u_\alpha$. The statement (2) can be verified by the direct
computation. \hfill$\Box$

\medskip

Define $v_\alpha^{(0)}:=1$ and 
$v_\alpha^{(k+1)}:=v_\alpha^{(k)}(\iota \bar{\iota})^k\iota(u_\alpha)$.
Then we have $v_\alpha^{(k)}\in \cN_{2k-1}$,  
$\alpha (\iota\bar{\iota})^k=\Ad (v^{(k)}_\alpha)\circ (\iota \bar{\iota})^k\alpha$, 
and $v_\alpha^{(k)*}$ satisfies 1-cocycle identity for $\alpha$.

Let $\alpha^{(k)}:=\Ad (v_\alpha^{(k)*})\circ \alpha$. Then $\alpha^{(k)}$
preserves Jones projections for a tunnel $\cM\supset \cN_1\supset \cdots \supset \cN_{2k}\supset
\cN_{2k+1}$. (Note Jones projections are given by
$\{(\iota\bar{\iota})^n(RR^*)\}_{n\geq 0}\cup
\{(\bar{\iota}\iota)^n(\bar{R}\bar{R}^*)\}_{n\geq 0}$.)
Hence $\alpha^{(k)}$ preserves this tunnel globally, and 
$\alpha^{(k)}(\cN_l'\cap \cM)=\cN_l'\cap \cM$ hold for all $0\leq l\leq 2k+1$.
We can see $\alpha^{(k)}|_{\cN_{2l+1}'\cap \cM}=
\alpha^{(l)}|_{\cN_{2l+1}'\cap \cM}$ for all $l\leq k$.

\begin{df}[\cite{Loi-auto}]
 The Loi invariant $\Phi(\alpha)$ of $\alpha$ is defined by 
$\Phi(\alpha)=\{\alpha^{(k)}|_{\cN_{2k+1}'\cap \cM}\}_{k}$.
\end{df}

By using the triviality of the Loi-invariant, we can construct a half
braiding unitary  $\mathcal{E}(\sigma,\alpha)\in
(\sigma\alpha,\alpha\sigma)$ for 
$[\sigma]\in \Delta$, 
$\alpha\in \mathrm{Ker}(\Phi)$. 
(The notion of a half braiding was introduced by Izumi \cite{Iz-LR},
inspired by the work of Xu \cite{Xu-newbraided}.)
Namely,  we have the following theorem.


\begin{thm}[{\cite[Theorem 2.1]{M-ExtIJM}}]\label{thm:halfbraiding}
 Let $\alpha\in \Aut(\cM,\cN)$. If $\Phi(\alpha)$ is trivial, then 
there exists a unitary $\mathcal{E}(\sigma,\alpha)$ for all 
$[\sigma]\in \Delta$ 
such that \\
$(\mathrm{i})$ $\mathcal{E}(\sigma,\alpha)\in (\sigma\alpha,\alpha\sigma)$. \\
$(\mathrm{ii})$ $\mathcal{E}(\sigma_1,\alpha)\sigma_1\left(\mathcal{E}(\sigma_2,\alpha)\right)T=
\alpha(T)\mathcal{E}(\sigma_3,\alpha)$ for any $[\sigma_i]\in \Delta$
and $T\in (\sigma_3,\sigma_1\sigma_2)$. \\
$(\mathrm{iii})$
$\mathcal{E}(\sigma,\alpha\beta)=\alpha(\mathcal{E}(\sigma,\beta))\mathcal{E}(\sigma,\alpha)$, 
$\mathcal{E}(\sigma,\Ad v)=v\sigma(v^*)$ for $v\in U(\cN)$.
\end{thm}

The second condition is a braiding fusion equation (BFE), and the third
condition means that $\mathcal{E}(\sigma,\alpha)^*$ is a 1-cocycle for $\alpha$.

Here we only explain how to construct $\mathcal{E}(\sigma,\alpha)$, and
the outline of proof. See \cite{M-ExtIJM} for detail of proof.

Let $[\sigma]\in {}_\cM\Delta_\cM$. 
Fix $n$ and an isometry 
$T\in (\sigma, (\iota\bar{\iota})^n)$. Define 
$W_T=\alpha(T^*)v_\alpha^{(n)}T$. 
It is clear that $W_T\in (\sigma\alpha,\alpha\sigma)$. 
By using $\Phi(\alpha)=1$, we can show that the definition of $W_T$ does not
depend on the choice of $T\in(\sigma,(\iota\bar{\iota})^n)$, and $W_T$ is a unitary.

Next we show $W_T$ does not depend on $n$. 
Take any $\pi \prec \sigma\iota $, and an isometry $S\in
(\pi,\sigma\iota)$. Then 
$\tilde{S}=\sqrt{\frac{d(\sigma)d(\iota)}{d(\pi)}}S^*\sigma(R)\in
(\sigma,\pi \bar{\iota})$ is an isometry \cite[Proposition 2.2]{Iz-CuntsII}.
We can easily verify $TS\tilde{S}\in (\sigma,(\iota\bar{\iota})^{n+1})$.
Again by the triviality of $\Phi(\alpha)$, we can show $W_T=W_{TS\tilde{S}}$. 

Combining these, we know that $W_T$ does not depend on $n$ and $T$. Hence 
$\mathcal{E}(\sigma,\alpha):=W_T$ is well-defined.
The condition (iii) follows from Lemma \ref{rem:1cocycle}.

To show (ii),  take $n,m$ and isometries $S_1\in (\sigma_1,(\iota\bar{\iota})^n)$,
$S_2\in (\sigma_2,(\iota\bar{\iota})^m)$. 
Then $S_3:=S_1\sigma_1(S_2)T\in (\sigma_3,(\iota\bar{\iota})^{n+m})$ is an isometry.
Then 
we can show $W_{S_1}\sigma_1(W_{S_2})T=\alpha(T)W_{S_3}$ by using
$\Phi(\alpha)
=1$. 
In a similar way, we can construct a half braiding
$\mathcal{E}(\sigma,\alpha)$ for each 
$[\sigma]\in {}_\cA\Delta_\cB$, $\cA,\cB\in\{\cN,\cM\}$. 
\hfill$\Box$

\smallskip 

\noindent
\textbf{Remark.}
We can extend $\mathcal{E}(\sigma, \alpha)$ for a reducible $\sigma$ as follows.
Let 
$\sigma=\sum_{i=1}^nw_i\sigma_i(x)w_i^*$ and  set 
$\mathcal{E}(\sigma,
\alpha):=\sum_{i}\alpha(w_i)\mathcal{E}(\sigma_i,\alpha)w_i^*$. Then 
BFE implies 
$\mathcal{E}(\rho\sigma,\alpha)=\mathcal{E}(\rho,\alpha)\rho(\mathcal{E}(\sigma,\alpha)
)$.

\section{Vanishing of obstruction in the orbifold construction}\label{sec:flat}

In this section, we make the following assumption. \\

\noindent
\textbf{Assumption.} \\
(A1) there exists $\alpha\in \Aut(\cM,\cN)$ such that 
$[\alpha]\in {}_\cM\Delta_\cM$, and $\alpha$ gives an outer action of
$\mathbb{Z}/n\mathbb{Z}$. \\
(A2) the Loi-invariant of $\alpha$ is trivial. \\
(A3) there exists a self-conjugate $[\rho]\in {}_\cM\Delta_\cM$ such that $[\alpha][\rho]=[\rho]$
and $[\rho]^2\succ [\rho]$. \\

As in \cite[Example 3.2]{Iz-CuntsI}, we can choose representatives of
$[\alpha]$ and $[\rho]$ such that 
$\alpha\rho=\rho$. In what follows, we fix this choice.

The crossed product inclusion $\cN\rtimes_\alpha
\mathbb{Z}/n\mathbb{Z}\subset 
\cM\rtimes_\alpha \mathbb{Z}/n\mathbb{Z}$ is called an orbifold
subfactor for $\cN\subset \cM$, and this construction is called the 
orbifold construction \cite{Kw-orbi}.

As seen in the previous section, 
 we have a half
braiding $\mathcal{E}(\sigma,\alpha)\in (\sigma\alpha,\alpha\sigma)$,
 $[\sigma]\in \Delta$. 
Once we get a half braiding, we can define an extension 
\[ 
\sigma\in \Mor(\cA,\cB)\rightarrow 
\tilde{\sigma}\in \Mor(\cA\rtimes_\alpha \mathbb{Z}/n\mathbb{Z},
\cB\rtimes_\alpha \mathbb{Z}/n\mathbb{Z})
\]
by setting 
$$\tilde{\sigma}(\lambda)=\mathcal{E}(\sigma, \alpha)^*\lambda,$$
where  $\lambda$ is an implementing unitary for $\alpha$, and
$\cA,\cB\in \{\cN,\cM\}$.
The condition (i) and (iii) in Theorem  \ref{thm:halfbraiding} imply 
that $\tilde{\sigma}$ indeed gives a
morphism. The condition (ii) implies that 
$(\sigma_3,\sigma_1\sigma_2)\subset 
(\tilde{\sigma}_3,\tilde{\sigma}_1\tilde{\sigma}_2)$. Thus the extension 
$\sigma\rightarrow \tilde{\sigma}$ preserves sector operation, and it 
is a special case of $\alpha$-induction studied in 
\cite{Xu-newbraided}, \cite{Boek-Evan-alphaI}, \cite{Boek-Evan-alphaII},
\cite{Boek-Evan-alphaIII}, \cite{BEK-alpha}. 
It is easy to see
$\hat{\alpha}\tilde{\sigma}=\tilde{\sigma}\hat{\alpha}$, 
where $\hat{\alpha}$ is the dual action given by
$\hat{\alpha}(\lambda)=\omega \lambda $, $\omega=e^{\frac{2\pi i}{n}}$.

It is trivial that
$\mathcal{E}(\alpha,\alpha)$ a scalar with
$\mathcal{E}(\alpha,\alpha)^n=1$,  but it may be non-trivial. We say 
$\mathcal{E}(\alpha,\alpha)$ an obstruction in the orbifold construction. This notion
comes from the following theorem.
\begin{thm}\label{thm:orbi}
Assume $\mathcal{E}(\alpha,\alpha)=1$. Then we have the following. \\
$(1)$ We have $\tilde{\alpha}=\Ad \lambda$. Thus
 $[\tilde{\sigma}]=[\widetilde{\alpha\sigma}]$ as sectors.
If $[\sigma]\ne
 [\alpha^k\sigma]$ for all $k=1,2,\cdots n-1$, then $\tilde{\sigma}$ is irreducible.
\\
$(2)$ 
We have
$(\tilde{\rho},\tilde{\rho})=
\{\sum_{k=0}^{n-1} a_k\lambda^k\mid
 a_k\in \mathbb{C}\}\cong \ell^\infty(\mathbb{Z}/n\mathbb{Z})$.
Therefore we have an irreducible decomposition $[\tilde{\rho}]=\oplus_{k=0}^{n-1}[\pi_i]$. 
Here $\pi_k$ is an irreducible sector corresponding to a minimal
 projection $p_k=n^{-1}\sum_{l=0}^{n-1}\omega^{kl}\lambda^l$. \\
$(3)$ Let $\hat{\alpha}$ be the dual action. Then
 $[\hat{\alpha}][\pi_k][\hat{\alpha}^{-1}]
=[\pi_{k+1}]$. Thus $d(\pi_k)=d(\rho)/n$.\\
$(4)$ If $n$ is odd, then all $[\pi_k]$ are self conjugate. 
If $n$ is even, then  $\overline{[\pi_k]}=[\pi_k]$, or 
$\overline{[\pi_k]}=\left[\pi_{k+\frac{n}{2}}\right]$ hold.
\end{thm}
\textbf{Proof}. 
(1) For $a\in \cM$, $\tilde{\alpha}(a)=\alpha(a)=\Ad (\lambda)(a)$ holds. For
$\lambda$,
$\tilde{\alpha}(\lambda)=\mathcal{E}(\alpha,\alpha)^*\lambda=\lambda
=\Ad (\lambda)(\lambda)$ by the assumption. Hence we have
$\tilde{\alpha}=\Ad (\lambda)$.

We show the latter statement for  $[\sigma]\in {}_\cM\Delta_{\cM}$. (Other cases can be verified in
the same way.)
Take $a=\sum_{k=0}^{n-1}a_k\lambda^k \in (\tilde{\sigma},\tilde{\sigma})$. 
For $x\in \cM$, 
\[
\sum_{k=0}^{n-1}\sigma(x)a_k\lambda^k=  \tilde{\sigma}(x)a=
 a\tilde{\sigma}(x)=\sum_{k=0}^{n-1}a_k\lambda^k\sigma(x)=
\sum_{k=0}^{n-1}a_k\alpha^k(\sigma(x))\lambda^k.
\]
Thus $a_k\in (\alpha^k\sigma,\sigma)$. By the assumption, $a_k=0$
for $1\leq k\leq n-1$, and $a_0\in \mathbb{C}1$. Hence $\tilde{\sigma}$
is irreducible. (For proof of this fact, the condition
$\mathcal{E}(\alpha,\alpha)=1$ is unnecessary.)\\
(2) Let $a=\sum_{k=0}^{n-1}a_k\lambda^k\in (\tilde{\rho},\tilde{\rho})$.
In a similar way as above, we get
$a_k\in(\alpha^k\rho,\rho)$. Since we have chosen $\alpha$ and $\rho$ so
that $\alpha\rho=\rho$, $a_k\in \mathbb{C}1$. 

If we apply BFE for $\sigma_1=\alpha$, $\sigma_2=\sigma_3=\rho$
 and $T=1\in(\rho,\alpha\rho)=(\rho,\rho)$, 
we get
$\alpha(\mathcal{E}(\rho,\alpha))=\mathcal{E}(\alpha,\alpha)^*\mathcal{E}(\rho,\alpha)$.  
Thus we have 
$\alpha(\mathcal{E}(\rho,\alpha))=\mathcal{E}(\rho,\alpha)$ by
assumption. 
Then we have
\[
 \tilde{\rho}(\lambda)a=\sum_{k=0}^{n-1}\mathcal{E}(\rho,\alpha)^*a_k\lambda^{k+1},
\]
and 
\[
 a\tilde{\rho}(\lambda)=\left(\sum_{k=0}^{n-1}a_k\lambda^k\right)\mathcal{E}(\rho,\alpha)^*\lambda
=\sum_{k=0}^{n-1}a_k\alpha^k(\mathcal{E}(\rho,\alpha)^*)\lambda^{k+1}=
\sum_{k=0}^{n-1}\mathcal{E}(\rho,\alpha)^*a_k\lambda^{k+1}.
\]
Thus $a\in (\tilde{\rho},\tilde{\rho})$, and we obtain
$(\tilde{\rho},\tilde{\rho})=\{\sum_{k=0}^{n-1} a_k\lambda^k\mid  a_k\in
\mathbb{C}\}$ . \\
(3) 
Fix an isometry $v_k$ with $v_kv_k^*=p_k$, and set $u=\hat{\alpha}(v_k^*)v_{k+1}$. 
Since $\hat{\alpha}(p_k)=p_{k+1}$, we can easily to see $u$ is a unitary.
Again by $\hat{\alpha}(p_k)=p_{k+1}$, we have
\begin{align*}
 \hat{\alpha}\pi_k\hat{\alpha}^{-1}(x)&=
\hat{\alpha}(v_k^*)
\hat{\alpha}\tilde{\rho}\hat{\alpha}^{-1}(x)\hat{\alpha}(v_k)=
\hat{\alpha}(v_k^*)\tilde{\rho}(x)\hat{\alpha}(v_k)=
\hat{\alpha}(v_k^*p_{k})\tilde{\rho}(x)\hat{\alpha}(p_kv_k)\\
&=\hat{\alpha}(v_k^*)p_{k+1}\tilde{\rho}(x)p_{k+1}\hat{\alpha}(v_k)
=uv_{k+1}^*\tilde{\rho}(x)v_{k+1}u  
=\Ad (u)\pi_{k+1}(x).
\end{align*}
(4) 
Assume $\overline{[{\pi}_0]}=[\pi_k]$ for some $k$. By considering the conjugate of 
$[\pi_k]=\hat{\alpha}^k[\pi_0]\hat{\alpha}^{-k}$, we get 
\[
 [\pi_0]=\overline{[\pi_k]}=\hat{\alpha}^k\overline{[\pi_0]}\hat{\alpha}^{-k}=\hat{\alpha}^k[\pi_k]\hat{\alpha}^{-k}
\]
and hence  $[\pi_0]=[\pi_{2k}]$ holds. Thus $k=0 $ if $n$ is odd. In this
case, all $[\pi_i]$ are self conjugate.
If $n=2n'$, then   $k=0$ or $k=n'$. The rest of statement can
be easily seen.
\hfill$\Box$

\medskip

\noindent
\textbf{Remark.} (1) It is easy to see $\tilde{\iota}$ is an inclusion
map $\cN\rtimes_\alpha\mathbb{Z}/n\mathbb{Z}\hookrightarrow
\cM\rtimes_\alpha\mathbb{Z}/n\mathbb{Z}$, and
$\widetilde{\iota\bar{\iota}}$ is a canonical endomorphism for this
inclusion. \\
(2) Even if we do not assume
$\mathcal{E}(\alpha,\alpha)=1$, we get similar results. For example,
the statement corresponding to (2) is the following: \\
Let $\mathcal{E}(\alpha,\alpha)$ be a primitive $l$-th root, and set $n/l=m$.
Then the irreducible decomposition of $\tilde{\rho}$ is
$\tilde{\rho}=\bigoplus_{k=0}^{m-1}[\pi_{k}]$, 

\medskip

By Theorem \ref{thm:orbi}, the graph change occurs by the orbifold construction under the
assumption $\mathcal{E}(\alpha,\alpha)=1$. (Also see the following examples.)

We would like to  determine when $\mathcal{E}(\alpha,\alpha)=1$. 
We have the following sufficient condition. 
\begin{thm}\label{thm:obst}
 Assume $m:=\dim(\rho,\rho^2)$ is relatively prime with $n$. Then
 $\mathcal{E}(\alpha,\alpha)=1$. 
\end{thm}
\textbf{Proof.} 
This proof is inspired by the computation in \cite{Iz-CuntsI}.
As explained in the proof of Theorem \ref{thm:orbi}, 
we have
$\alpha(\mathcal{E}(\rho,\alpha))=\mathcal{E}(\alpha,\alpha)^*\mathcal{E}(\rho,\alpha)$. 
Put $U=\mathcal{E}(\rho,\alpha)^*$, and  
$a=\mathcal{E}(\alpha,\alpha)$. Of course we have 
$a^n=1$.

Since $(\rho\alpha,\alpha\rho)=(\rho\alpha,\rho)$, we have
$\rho\alpha=\Ad U\rho$, and hence $U\in(\rho,\rho\alpha)$.

We also have $U\in (\rho^2,\rho^2)$ due to $\alpha\rho=\rho$. 
Let $z\in (\rho^2,\rho^2)$ be a minimal central projection corresponding
to the irreducible component $\rho$ of $\rho^2$.
Fix an orthonormal basis $\{T_i\}_{i=1}^m\subset (\rho,\rho^2)$. 
We have $\sum_{i=1}^mT_iT_i^*=z$, and 
\[
U=
\sum_{i=1}^m d_{ij}T_iT_j^*+(1-z)U=
(T_1,\cdots , T_m)D
\left(\begin{array}{c}
 T_1^*\\ \vdots  \\ T_m^*
\end{array}\right)
+(1-z)U
\]
for some unitary matrix $D=(d_{ij})\in M_m(\mathbb{C})$.

Since $\alpha$ acts on $(\rho,\rho^2)$ as a unitary and $\alpha(z)=z$, there exists a 
unitary matrix $V\in M_m(\mathbb{C})$ such that 
$$(\alpha(T_1),\cdots \alpha(T_m))=(T_1,\cdots, T_m)V.$$

By the condition $\alpha(U)=aU$, 
we get $VDV^*=aD$. By taking determinant, 
we get $\mathrm{det}(D)=a^m\mathrm{det}(D)$. Hence we get $a^m=1$. Since
$m$ and $n$ are relatively prime, we have $a=1$. \hfill$\Box$

\medskip
\noindent
\textbf{Remark.} (1) 
In the above proof, if we have $\mathrm{Tr}(D)\ne 0$,
then $a=1$ follows immediately. However it does not seem easy to
determine $\mathrm{Tr}(D)\ne 0$ in general case.  \\
(2) Let us assume that  ${}_\cM\Delta_\cM$ is a near group category, i.e.,
${}_\cM\Delta_\cM=\{[\alpha_g]\}_{g\in G}\cup\{[\rho]\}$ for some action $\alpha$
of a finite group $G$. 
Gannon and Evans show that $\dim(\rho,\rho^2)$ is $|G|-1$ or multiple of
$|G|$ in \cite{EvGan-neargrp}. (In the former case, $|G|$ is prime.)   

\medskip

\begin{exam}\label{exam:Aodd}
\upshape
Let $\cN\subset \cM$ be a Jones subfactor with 
principal graph $A_{4n-3}$ \cite{J-ind}. 
Sectors of $\cN\subset \cM$ appear as follows:
\[
 [\rho_0]-[\rho_1]-\cdots -[\rho_{2n-2}]-\cdots -[\rho_{4n-4}]
\]
Here $\rho_0=\id_{\cM}$, and $\rho_1=\iota$. Then $\alpha:=\rho_{4n-4}$ 
is an automorphism of $\cN\subset \cM$
with period 2. 
It is well known that all $[\rho_{2k}]$ are self conjugate, and 
$[\alpha][\rho_k]=[\rho_{4n-4-k}]$ hold for $0\leq k\leq 2n-2$. (See \cite{Iz-fusion}.)
In particular, the sector $[\rho_{2n-2}]$ is self conjugate,  and satisfies 
$[\alpha][\rho_{2n-2}]=[\rho_{2n-2}]$. Moreover, we have $\dim(\rho_{2n-2},\rho_{2n-2}^2)=1$. 
Thus the assumption of Theorem \ref{thm:obst} is 
satisfied. 
By Theorem \ref{thm:orbi}, we have $[\tilde{\rho}_k]=[\tilde{\rho}_{4n-4-k}]$
 for $0\leq k< 2n-2$, and $[\tilde{\rho}_{2n-2}]=[\pi_0]\oplus [\pi_1]$. 
Hence the principal graph of the orbifold subfactor is  a
Dynkin diagram $D_{2n}$.  
\end{exam}

\begin{exam}\label{exam:ExtE6}
\upshape
Let $\cN\subset \cM$ be a  subfactor with a
principal graph $E^{(1)}_6$. Then 
${}_\cM\Delta_\cM$ is $\{[\id],[\alpha],[\alpha]^2,[\rho]\}$ with the following fusion rule.
\[
[\alpha]^3=[\id],\, [\alpha][\rho]=[\rho][\alpha]=[\rho],\,
 [\rho]^2=[\id]\oplus [\alpha]\oplus [\alpha^2]\oplus 2[\rho]
\]
We can take $\alpha$ as $\alpha^3=\id$, and $\alpha$ gives an outer
action of $\mathbb{Z}\slash 3\mathbb{Z}$ on $\cN\subset \cM$ with trivial Loi invariant.
Thus Theorem \ref{thm:obst} can be applied. 
In this case, we have
\[
[\tilde{\id}]=[\tilde{\alpha}]=[\tilde{\alpha}^2],\,\,\, 
[\tilde{\rho}]=[\pi_0]\oplus [\pi_1]\oplus[\pi_2]. 
\]
Therefore, the principal graph of an orbifold subfactor is
$D^{(1)}_{4}$. (Note the statistical dimension of $\rho$ is $d(\rho)=3$.)
There exist two subfactors with principal graph 
$D^{(1)}_{4}$, which arise as the crossed product by $\mathbb{Z}/4\mathbb{Z}$, and 
$\mathbb{Z}/2\mathbb{Z}\times \mathbb{Z}/2\mathbb{Z}$, respectively.
The condition (4) of Theorem \ref{thm:orbi} implies that the orbifold
subfactor is the crossed product by $\mathbb{Z}/2\mathbb{Z}\times \mathbb{Z}/2\mathbb{Z}$.

The subfactor treated  above 
is an $SU(3)_{3}$
 subfactor. (See \cite{Wen-Hecke}, \cite{EK-HeckeI} for $SU(N)_{l}$ subfactors.)
We can apply our main theorem to $SU(3)_{3k}$-subfactors
for $k+1\not\equiv 0(\mathrm{mod}\,\, 3)$. 
Indeed, there exist $\alpha\in \Aut(\cM,\cN)$ with
 $\alpha^3=\id$, $\Phi(\alpha)=\id$, and a unique self conjugate sector $[\rho]$ fixed by $\alpha$
 for an $SU(3)_{3k}$ subfactor $\cN\subset \cM$. (The sector $[\rho]$ corresponds to a
 young diagram $(2k,k,0)$.)
We have $m=k+1$ for an $SU(3)_{3k}$-subfactor by
applying  the Littlewood-Richardson rule for $SU_q(3)_{3k}$ \cite{GoodWen-LittRechard}.
(When $k+1\equiv 0(\mathrm{mod} \,\,3 )$, we can
not apply Theorem \ref{thm:obst}. However, obstruction vanishes in this
 case \cite{EK-HeckeI}.) 
\end{exam}

If the assumption of Theorem \ref{thm:obst} is not satisfied, 
$\mathcal{E}(\alpha,\alpha)$ may be fail to be  1. 

\begin{exam}[{\cite[Example 3.4]{Iz-CuntsI}}]\label{exam:E6}
\upshape
Let $\cN\subset \cM$ be a subfactor with principal graph $E_6$. 
Then ${}_\cM\Delta_{\cM}=\{[\id],[\alpha], [\rho]\}$ and 
they obey the following fusion rule:
\[
  [\rho^2]=[\id]\oplus [\alpha] \oplus 2[\rho], \,\, [\alpha]^2=[\id], \,\,
[\alpha][\rho]=[\rho].
\]
We have $\Ker(\Phi)=\Aut(\cM,\cN)$, and hence the Loi invariant of
 $\alpha$ is trivial. 
Izumi showed that one can take  $S_1\in (\id,\rho^2)$, $S_2\in (\alpha,\rho^2)$,
$S_3,S_4\in (\rho,\rho^2)$ and $U\in (\rho,\rho\alpha)$ as 
\[
 S_2=\alpha(S_1), \, \alpha(S_3)=S_3,\, \alpha(S_4)=-S_4,\,
U=S_1S_1^*-S_2S_2^*+S_3S_4^*+S_4S_3^*.
\]
Then $\alpha(U)=-U$, and
hence $\mathcal{E}(\alpha,\alpha)=-1$. In this case, the graph change
does not occur by the orbifold construction. 
\end{exam}

\ifx\undefined\bysame
\newcommand{\bysame}{\leavevmode\hbox to3em{\hrulefill}\,}
\fi

\end{document}